\documentclass[10pt]{article}
\usepackage{amsmath,amsfonts}
\usepackage{verbatim}
\usepackage{relsize}
\normalbaselines
\newtheorem{teor}{Theorem}
\newtheorem{prop}[teor]{Proposition}

\newtheorem{coro}[teor]{Corollary}
\newtheorem{rem}[teor]{Remark}

\newenvironment{demo}{\rm \trivlist \item[\hskip \labelsep{\it
      Proof}.]}{\nopagebreak \hfill $\square$ \endtrivlist}

\title{On uniqueness of the foliation by comoving observers restspaces of a Generalized Robertson Walker spacetime}

\author{Jos\'e A. S. Pelegr\'in${}^{a}$, Alfonso
Romero${}^{a}$ and Rafael M. Rubio${}^{b}$ \\[6mm]
${}^a$ Departamento de Geometr\'\i a y Topolog\'\i a, \\ [0.5mm]
Universidad de Granada, 18071 Granada, Spain \\ E-mails\textup{:
\texttt{jpelegrin@ugr.es, aromero@ugr.es}} \\[3mm]
${}^b$ Departamento de Matem\'aticas, Campus de Rabanales, \\[0.5mm] Universidad de
C\'ordoba, 14071 C\'ordoba, Spain,\\[0.5mm] E-mail\textup{: \texttt{rmrubio@uco.es}}\\[3mm]}

\date{}

\begin{document}

\maketitle

\thispagestyle{empty}

\begin{abstract}
A characterization of the foliation by spacelike slices of an $(n+1)$-dimensional spatially closed Generalized Robertson-Walker spacetime is given by means of studying 
a natural mean curvature type equation on spacelike graphs. Under some natural assumptions,
of physical or geometric nature, all the entire solutions of such an equation are
obtained. In particular, the case of entire spacelike graphs in de Sitter spacetime is faced and completely solved by means of a new application
of a known integral formula.
\end{abstract}
\vspace*{5mm}




\section{Introduction}

The study of spacelike hypersurfaces has a long and fruitful history in General Relativity \cite{MT}. Indeed, the existence of constant mean curvature spacelike (and in particular maximal) hypersurfaces is necessary for the study of the structure of singularities in the space of solutions
of Einstein's equations \cite{AMM}. Also, the deep understanding of this kind of hypersurfaces is essential to prove the positivity of the gravitational mass \cite{SY}.
They are also interesting for numerical Relativity, where maximal hypersurfaces are used for integrating forward in time \cite{JKG}.

Among classical papers dealing with uniqueness results for constant mean curvature (CMC) spacelike hypersurfaces we can point out
\cite{BF}, \cite{Ch} and \cite{MT}, although a previous relevant
result in this direction was the proof of the Calabi-Bernstein
conjecture \cite{Ca} for the $(n+1)$-dimensional Lorentz-Minkowski
spacetime given by Cheng and Yau \cite{CY}. More recently, these kind of uniqueness results have been extended to a wide variety of spacetimes (see for instance \cite{A-R-S1}, \cite{PRR1} and \cite{PRR2}).

In this article, we will study spacelike hypersurfaces in the family of cosmological models known as Generalized Robertson-Walker (GRW) spacetimes. GRW spacetimes are warped products whose negative definite base represents a universal time and whose fiber is an arbitrary Riemannian manifold (see Section 2). These spacetimes play a relevant role in Physics as well as extend the classical notion of Robertson-Walker spacetime to the case when the fiber does not necessarily have constant sectional curvature \cite{A-R-S1}. On the other hand, deformations of the metric on the fiber of classical Robertson-Walker spacetimes as well as GRW spacetimes with a time-dependent conformal change of metric also fit into the class of GRW spacetimes.  In fact, in this paper we will deal with spatially closed GRW spacetimes, that is, GRW spacetimes whose fiber is compact \cite[Prop. 3.2]{A-R-S1}.

In any GRW spacetime there is a remarkable family of observers, the so-called comoving observers, which are collected as the integral curves of the reference frame $\partial_t$ \cite{SW}. The volume of the physical space measured by any comoving observer $\gamma$ at an instant $t$ of its proper time is given by

\begin{equation}
\label{volsli}
\textrm{Vol}_{\gamma}(t) = f(t)^n \textrm{Vol}(F),
\end{equation}

\noindent where $\textrm{Vol}(F)$ is the volume of the fiber $F$ and $f$ represents the warping function. Indeed, the comoving observers determine a foliation of the GRW spacetime by spacelike slices, which are their restspaces. From (\ref{volsli}), we see that the volume change that the comoving observers measure for the spacelike slices is given by the function $n f(t)^{n-1}f'(t)$. In general, we can consider the volume variation of the leaves of a foliation of the spacetime by closed spacelike hypersurfaces and ask ourselves the following question:

\begin{quote}
\begin{center}
\textit{Does there exist another (local) foliation of a spatially closed GRW spacetime\\
 by closed spacelike hypersurfaces with the same volume variation\\
 that the comoving observers measure for their restspaces?}
\end{center}
\end{quote}

Since a compact spacelike hypersurface in a spatially closed GRW spacetime with simply connected fiber is a graph \cite[Prop. 3.3]{A-R-S1}, to obtain the volume variation of the leaves of a foliation by spacelike hypersurfaces we will study the volume variation of spacelike graphs. Thus, let us consider in $\overline{M}$ the graph $\Sigma_u = \{(u(p), p) : p \in F \}$, where $F$ is the fiber, $u \in C^\infty(F)$ and $u(F)\subset I$. The induced metric on $F$ from the Lorentzian metric on $\overline{M}$ (see Section 2), via the graph $\Sigma_u$ is given by

\begin{equation}
\label{megra}
 g_u = -du^2 + f(u)^2 g,
\end{equation}

\noindent where $g$ is the metric on $F$. The metric $g_u$ is positive definite (i.e., $\Sigma_u$ is spacelike) if and only if $u$ satisfies $|Du| < f(u)$. In this case

\begin{equation}
\label{normalvec}
N = \frac{1}{f(u) \sqrt{f(u)^2 - |Du|^2}} \left( f(u)^2 , Du \right)
\end{equation}

\noindent is a future pointing unit normal vector field on $\Sigma_u$. Hence, the volume of this spacelike graph is

\begin{equation}
\label{volgra}
\textrm{Vol}(\Sigma_u) = \int_F f(u)^{n-1} \sqrt{f(u)^2 - |Du|^2} dV^F,
\end{equation}

\noindent where $dV^F$ is the canonical measure on $F$ given by $g$.

Taking (\ref{volsli}) and (\ref{volgra}) into account, we can easily see that the spacelike graphs that answer our question will be the ones defined by a function that solves the following variational problem:

Let $f: I \longrightarrow \mathbb{R}$ be a positive smooth function on an open interval $I \subseteq \mathbb{R}$ and let $(F, g)$ be an $n$-dimensional compact Riemannian manifold. Consider the class of smooth real valued functions $u$ on $F$ such that $u(F)\subset I$ and the length of the gradient of $u$ in $(F,g)$ satisfies $|Du|<f(u)$. On this class consider the functional

\begin{equation}
\label{functional}
\mathcal{I}(u) := \int_F \left[ f(u)^{n-1} \sqrt{f(u)^2 - |Du|^2} - f(u)^n \right] dV^F.
\end{equation} 

Note that the stationary points of this functional will determine foliations of $\overline{M}$ by spacelike graphs with the same volume variation than the foliation by spacelike slices. For stationary points of this functional, the Euler-Lagrange equation is written as

\begin{samepage}

\begin{equation}
\label{eul1} \tag{E.1} 
H(u) = \frac{f'(u)}{f(u)}, 
\end{equation}
\begin{equation}
\label{eul2} \tag{E.2}
|Du|<f(u),
\end{equation}

\end{samepage}

\noindent being $H(u)$ the smooth function given by

\begin{equation}
\label{curvmed}
H(u) =  {\rm div} \left( \frac{Du}{n f(u) \sqrt{f(u)^2 - |Du|^2}} \right) 
 +\frac{f'(u)}{n \sqrt{f(u)^2 - |Du|^2}} \left( n + \frac{|Du|^2}{f(u)^2} \right),
\end{equation}

\noindent where ${\rm div}$ represents the divergence operator in $(F, g)$.

Indeed, $H(u)$ represents the mean curvature of $\Sigma_u$ in $\overline{M}$ with respect to $N$. Thus, the variational problem that answers our question is equivalent to a mean curvature prescription problem. 

The study of the prescribed mean curvature Dirichlet problem began in the Euclidean context with Serrin's work \cite{Se}, who found necessary and sufficient conditions for its solvability. In Lorentz-Minkowski spacetime, a crucial result on necessary and sufficient conditions for the existence of smooth strictly spacelike solutions of the Dirichlet problem for the mean curvature operator was given in \cite{BS} by Bartnik and Simon. Gerhardt \cite{Ge1} also studied the existence of closed spacelike hypersurfaces with prescribed mean curvature in a globally hyperbolic Lorentzian manifold with a compact Cauchy hypersurface. Later, Ecker and Huisken gave a simpler proof in \cite{EH} for some of the results in \cite{Ge1} using an evolutionary approach. However, they needed to impose restrictive assumptions to obtain their results. Finally, Gerdhart showed in \cite{Ge2} that the evolutionary method could be used to prove his earlier result without imposing additional assumptions. More recently, several authors have also studied mean curvature prescription problems in Lorentz-Minkowski spacetime (see \cite{Azz}, \cite{BJM}, \cite{CCR} and references therein) as well as in more general ambiences \cite{dRT}. Note that these papers deal with local solutions, whereas we will give uniqueness results for global solutions of our problem.

Hence, in this article we will prove some uniqueness results for entire solutions of a prescribed mean curvature problem in a more general family of spacetimes, namely, Generalized Robertson-Walker spacetimes with compact fiber. In fact, the prescription function will be the Hubble function $f'/f$, which has an important physical meaning in General Relativity. Therefore, we will study the problem of determining the compact spacelike hypersurfaces immersed in this class of spacetimes whose mean curvature function at each point is the same than the mean curvature of the spacelike slice that passes through that point. As we have previously seen, this problem is equivalent to find the observers that measure a physical space with the same volume variation that the comoving ones measure for the spacelike slices when the fiber is simply connected.

We can make some comments on the meaning of (E): \textbf{1)} This equation represents a mean curvature prescription problem, i.e., we are dealing with a nonlinear elliptic problem. \textbf{2)} The right member of (\ref{eul1}) is, at each point $(u(p), p) \in \overline{M}$, the mean curvature of the spacelike slice $t=u(p)$. \textbf{3)} However, this is not a comparison assumption between extrinsic quantities of two spacelike hypersurfaces in $\overline{M}$, since the right member of (\ref{eul1}) corresponds to a spacelike slice that changes at each point. \textbf{4)} When the warping function is constant, (E) turns into the well-known maximal hypersurface equation.

The existence of local solutions of the nonlinear elliptic problem (E) is ensured by \cite[Thm. 4.3]{Kaz}. Nevertheless, here we are interested in global (i.e., entire) solutions of the problem. Concretely, in this article we will state several uniqueness results for entire solutions of (E). 

The main achievement in this article is the following uniqueness result for equation (E)

\begin{quote}
\textit{Let $F$ be an $n$-dimensional compact Riemannian manifold and let $f: I \longrightarrow ]0, \infty[$ be a smooth function such that $f'$ is signed. Then, the only entire solutions of equation} (E) \textit{are the constant functions.}
\end{quote}

This result is in fact a consequence of Corollary \ref{calbeleq}. Indeed, problems related to the one solved in Corollary \ref{calbeleq} have been previously studied in the particular case where the fiber is a complete noncompact Riemannian surface (see \cite{CRR}, \cite{RR} and \cite{RRS3}) as well as for compact 2-dimensional Riemannian fiber in \cite{RR2}. The techniques used in these cases strongly depended on the rich conformal geometry in the 2-dimensional case. Here, we deal with this problem when the fiber $F$ is a compact Riemannian manifold of arbitrary dimension. Our approach to equation (E) is purely geometric, since we study first a parametric version of the problem (Theorem \ref{teoleq}) that will enable us to solve (E).

Furthermore, we will also solve equation (E) in the case where $F=\mathbb{S}^n$, the unit round n-sphere and $f: \mathbb{R} \longrightarrow \mathbb{R}$, $f(t) = \cosh(t)$, i.e., for entire spacelike graphs in de Sitter spacetime $\overline{M} = \mathbb{R} \times_{\cosh(t)} \mathbb{S}^n$. In order to do so, we will deal with the more general problem of deciding when a compact spacelike hypersurface $M$ in an Einstein GRW spacetime whose mean curvature function $H$ satisfies $H = \frac{f'(\tau)}{f(\tau)}$, being the right member the restriction of the Hubble function to $M$, is a spacelike slice (Theorem \ref{teointeins}). We have done this by means of a new application of a well-known integral formula (see Section \ref{seceins}). Thus, we have (Corollary \ref{nonpards})

\begin{quote}
\textit{The only entire solutions on $\mathbb{S}^n$ of}
$$ H(u) = \tanh(u),$$
$$ |Du| < \cosh(u)$$
\noindent \textit{are the constant functions.} 
\end{quote}
 
Finally, in Section \ref{last} we will deal with GRW spacetimes obeying the so-called Null Convergence Condition (NCC). This curvature assumption is satisfied by a great deal of relevant cosmological models, including some describing inflationary scenarios. Moreover, it appears as a natural assumption when studying spacetimes' singularities. Following an analogous procedure as in the previous section, we classify the compact spacelike hypersurfaces satisfying $H = \frac{f'(\tau)}{f(\tau)}$ in certain GRW spacetimes obeying NCC (Theorem 
\ref{teointeins2}). As a consequence we get (Corollary \ref{ultimo})
\begin{quote}
\textit{
Let $F$ be an $n$-dimensional compact Riemannian manifold and let $f:I\longrightarrow ]0,\infty[$ be a smooth function such that $\log f$ is convex. If the Ricci curvature of $F$ is strictly bounded from below by $(n-1)\sup\, \big(f^2(\log f)''\big)$ or $\log f$ is strictly convex, then the only entire solutions of equation {\rm (E)} are the constant functions.}
\end{quote}

\section{Preliminaries}
\label{s2} 

Let $(F,g)$ be an $n(\geq 2)$-dimensional (connected) Riemannian manifold, $I$ an open interval in $\mathbb{R}$ and $f$ a positive smooth function defined on $I$. Now, consider the product manifold $\overline{M} = I \times F$ endowed with the Lorentzian metric

\begin{equation}
\label{metr}
\overline{g} = -\pi^*_{_I} (dt^2) +f(\pi_{_I})^2 \, \pi_{_F}^* (g), 
\end{equation}

\noindent where $\pi_{_I}$ and $\pi_{_F}$ denote the projections onto $I$ and
$F$, respectively. The Lorentzian manifold $(\overline{M}, \overline{g})$ is a warped product (in the sense of \cite{O'N}) with base $(I,-dt^2)$, fiber $(F,g)$ and warping function $f$. If we endow $(\overline{M}, \overline{g})$ with the time orientation induced by $\partial_t := \partial / \partial t$ we can call it an $(n+1)$-dimensional Generalized Robertson-Walker (GRW) spacetime.

The distinguished vector field $K: =~ f({\pi}_I)\,\partial_t$ is timelike and future pointing. From the relation between the
Levi-Civita connection of $\overline{M}$ and those of the base and
the fiber \cite[Cor. 7.35]{O'N}, it follows that

\begin{equation}\label{conexion} \overline{\nabla}_XK =
f'({\pi}_I)\,X
\end{equation}

\noindent for any $X\in \mathfrak{X}(\overline{M})$, where $\overline{\nabla}$
is the Levi-Civita connection of the Lorentzian metric
(\ref{metr}). Hence, $K$ is conformal and its metrically equivalent $1$-form is closed.

Given an $n$-dimensional manifold $M$, an immersion $\psi: M
\rightarrow \overline{M}$ is said to be spacelike if the
Lorentzian metric (\ref{metr}) induces, via $\psi$, a Riemannian
metric $g_{_M}$ on $M$. In this case, $M$ is called a spacelike
hypersurface. We will denote by $\tau:=\pi_I\circ \psi$ the
restriction of $\pi_I$ along $\psi$.

The time-orientation of $\overline{M}$ allows to take, for each
spacelike hypersurface $M$ in $\overline{M}$, a unique unitary
timelike vector field $N \in \mathfrak{X}^\bot(M)$ globally defined
on $M$ with the same time-orientation as $\partial_t$. Hence, from the wrong way Cauchy-Schwarz inequality, (see \cite[Prop. 5.30]{O'N}) we obtain $\overline{g}(N,\partial_t)\leq -1$ and $\overline{g}(N,\partial_t)= -1$ at a point $p\in M$ if and only if $N(p) = \partial_t(p)$. We will denote by $A$ the shape operator associated to $N$.
Then, the mean curvature function associated to $N$ is given
by $H:= -(1/n) \mathrm{trace}(A)$. As it is well-known, the mean
curvature is constant if and only if the spacelike hypersurface is,
locally, a critical point of the $n$-dimensional area functional for
compactly supported normal variations, under certain constraints of
the volume \cite{BC}. When the mean curvature vanishes identically, the
spacelike hypersurface is called a maximal hypersurface \cite{MT}.

In any GRW spacetime $\overline{M}= I \times_f F$ there
is a remarkable family of spacelike hypersurfaces, namely its
spacelike slices $\{t_{0}\}\times F$, $t_{0}\in I$. It can be easily
seen that a spacelike hypersurface in $\overline{M}$ is a (piece of)
spacelike slice if and only if the function $\tau$ is constant.
Furthermore, a spacelike hypersurface $M$ in $\overline{M}$ is a (piece
of) spacelike slice if and only if $M$ is orthogonal to $\partial_t$. The shape operator of the spacelike slice
$t=t_{0}$ is given by $A=-f'(t_{0})/f(t_{0})\,\mathbb{I}$, where $\mathbb{I}$
denotes the identity transformation, and therefore its (constant)
mean curvature is equal to the Hubble function at $t_{0}$, i.e., $H= f'(t_{0})/f(t_{0})$. Thus, a spacelike slice
is maximal if and only if $f'(t_{0})=0$ (and hence, totally
geodesic).

Let $\psi: M \rightarrow \overline{M}$ be an $n$-dimensional
spacelike hypersurface immersed in a Generalized Robertson-Walker spacetime
$\overline{M}= I \times_f F$. If we denote by
$\partial_t^T:= \partial_t+\overline{g}(N,\partial_t)N$
the tangential component of $\partial_t$ along $\psi$, then it is
easy to check that the gradient of $\tau:=\pi_I\circ \psi$ on $M$ is
\begin{equation}\label{part}
\nabla \tau=-\partial_t^T.
\end{equation}

Moreover, since the tangential
component of $K$ along $\psi$ is given by $K^T=K+\overline{g}(K,N)N$, a direct computation from
(\ref{conexion}) gives
\begin{equation}\label{gradcosh}
\nabla \overline{g}(K,N)=-AK^T.
\end{equation}

On the other hand, if we represent by $\nabla$ the Levi-Civita
connection of the metric $g_{_M}$, then the Gauss and Weingarten
formulas for the immersion $\psi$ are given, respectively, by
\begin{equation}\label{GF}
\overline{\nabla}_X Y=\nabla_X Y-g_{_M}(AX,Y)N
\end{equation}
and
\begin{equation}\label{WF}
AX=-\overline{\nabla}_X N,
\end{equation}
where $X,Y\in\mathfrak{X}({M})$. Then,  taking the tangential component in
(\ref{conexion}) and using (\ref{GF}) and (\ref{WF}), we get
\begin{equation}\label{KT}
\nabla_X K^T=-f(\tau)\overline{g}(N,\partial_t)AX+f'(\tau)X,
\end{equation}
where $X\in\mathfrak{X}({M})$ and $f'(\tau):=f'\circ \tau$. Since
$K^T=f(\tau)\partial_t^T$, it follows from (\ref{part}) and
(\ref{KT}) that the Laplacian of $\tau$ on $M$ is
\begin{equation}\label{laptau}
\Delta \tau=-\frac{f'(\tau)}{f(\tau)}\{n+|\nabla
\tau|^2\}-nH\overline{g}(N,\partial_t).
\end{equation}

\section{Uniqueness results when the warping function is monotone}

Now, we are able to obtain some uniqueness results. In fact, we may wonder if there are any other compact spacelike hypersurfaces apart from the spacelike slices satisfying $H = \frac{f'(\tau)}{f(\tau)}.$ Indeed, we will go one step further

\begin{teor}
\label{teoleq}
The only compact spacelike hypersurfaces in a GRW spacetime whose mean curvature function satisfies $H \leq \frac{f'(\tau)}{f(\tau)}$ with $f'(\tau) \leq 0$ (resp. $H \geq \frac{f'(\tau)}{f(\tau)}$ with $f'(\tau) \geq 0$) are the spacelike slices.
\end{teor}

\begin{demo}
Let us consider a primitive function $\mathcal{F}$ of $f$ and write $\mathcal{F}(\tau)$ for the restriction of $\mathcal{F} \circ \pi_I $ on M. Note that $\nabla \mathcal{F}(\tau) = f(\tau) \nabla \tau.$ Using (\ref{laptau}) we can compute this function's Laplacian

\begin{equation}
\label{lapf}
\Delta \mathcal{F}(\tau) = f(\tau) \Delta \tau + f'(\tau) |\nabla \tau|^2 = -n f'(\tau) - n f(\tau) H \overline{g}(N, \partial_t).
\end{equation}

From our assumptions and (\ref{lapf}) we have 

$$\Delta \mathcal{F}(\tau) \leq -n f'(\tau) (1 + \overline{g}(N, \partial_t)) \leq 0.$$

Hence, the compactness of $M$ implies that $\mathcal{F}$ must be constant. Consequently, $\nabla \mathcal{F}(\tau) = f(\tau) \nabla \tau = 0$ and $M$ is a spacelike slice. We can prove analogously the other statement.
\end{demo}

In particular, we have obtained 

\begin{quote}
\label{teoeq}
\textit{The only compact spacelike hypersurfaces in a GRW spacetime such that $f'(\tau)$ is signed and whose mean curvature function $H$ satisfies $H = \frac{f'(\tau)}{f(\tau)}$ are the spacelike slices.}
\end{quote}

\begin{rem}
\label{remarkexpand}
\normalfont
Following the ideas in \cite[Chap. 12]{O'N}, consider in the spacetime $\overline{M}= I \times_f F$ for $p \in F$ the integral curve $\gamma_p$ of the galactic flow $\partial_t$ through $p$, i.e., each $\gamma_p$ represents the evolution of a galaxy in $\overline{M}$. The time function $t$ may be seen as the common proper time of every $\gamma_p$. On the spacelike slice $t= t_0$ the distance with respect to the induced Riemannian metric between two galaxies at the same time $(t_0, p)$, $(t_0, q)$ is $f(t_0) d_F(p,q)$, where $d_F$ is the distance associated to $(F,g)$. Thus, the assumption $f' \geq 0$ (resp., $f' \leq 0$) means that any two galaxies are not coming together (resp., not spreading out).
\end{rem}

Therefore, the previous result gives 

\begin{quote}
\label{cornonexp}
\textit{The only compact spacelike hypersurfaces in a non-expanding (resp. non-contracting) GRW spacetime whose mean curvature function satisfies $H \leq \frac{f'(\tau)}{f(\tau)}$ (resp. $H \geq \frac{f'(\tau)}{f(\tau)}$) are the spacelike slices.}
\end{quote}

\begin{rem}
\label{remcompa}
\normalfont
Notice that in Theorem \ref{teoleq}: {\bf a)} There is no curvature assumption on the spacetime . In fact, if we wanted to obtain a similar result using the technique suggested in \cite[Rem. 2.3]{RR2}, we would need to impose these kind of curvature restrictions. {\bf b)} The compactness assumption on the hypersurface  cannot be weakened to completeness. Indeed, a trivial obstruction is given by any spacelike hyperplane in the $(n+1)$-dimensional Lorentz-Minkowski spacetime. Even more, there exist complete spacelike hypersurfaces in the $(n+1)$-dimensional steady state spacetime satisfying $H=1$ different from the spacelike slices \cite{AA}. {\bf c)} Despite $f'$ not being signed, we can add some assumption to the spacelike hypersurface in order to apply Theorem \ref{teoleq} (see Proposition \ref{prodesit}). {\bf d)} If the inequality satisfied by the mean curvature function is strict the result turns into a nonexistence one.
\end{rem}

As a consequence of Theorem \ref{teoleq}, we have the following uniqueness result

\begin{coro}
\label{calbeleq}

Let $F$ be an $n$-dimensional compact Riemannian manifold and let $f: I \longrightarrow ]0, \infty[$ be a smooth function such that $f' \leq 0$ (resp. $f' \geq 0$). Then, the only entire solutions of

$$ H(u) \leq \frac{f'(u)}{f(u)} \ \left( \text{resp.} \ H(u) \geq \frac{f'(u)}{f(u)} \right),$$
$$ |Du| < f(u)$$

\noindent are the constant functions. 

\end{coro}

The previous result contains the following nonparametric Calabi-Bernstein type result for (E).

\begin{quote}
\label{calaber1}
\textit{Let $F$ be an $n$-dimensional compact Riemannian manifold and let $f: I \longrightarrow ]0, \infty[$ be a smooth function such that $f'$ is signed. Then, the only entire solutions of equation} (E) \textit{are the constant functions.}
\end{quote}

We end this section highlighting that in Theorem \ref{teoleq} $H$ and $f'(\tau)$ are equally signed. However, if they have different sign we get

\begin{prop}
\label{teoHf}
The only compact spacelike hypersurfaces in a non-contracting (resp. non-expanding) GRW spacetime whose mean curvature function satisfies $H \leq 0$ (resp., $H \geq 0$) are the totally geodesic spacelike slices.
\end{prop}

\begin{demo}
From our assumptions we have $H f'(\tau) \leq 0$. Now, using (\ref{laptau}) we obtain that $\Delta \tau$ is signed. Since the hypersurface is compact, it must be a spacelike slice. Moreover, the only spacelike slices that satisfy this assumption are the totally geodesic ones.
\end{demo}

\begin{rem}
\label{rhs}
\normalfont
It should be noticed that the technique used to prove Theorem \ref{teoleq} may be extended to deal with other elliptic operators on the spacelike hypersurface $M$ appart from the Laplacian. In fact, associated to the $k$th Newton transformation on $M$, $P_k$, $0 \leq k \leq n$, there exists a second order linear differential operator $L_k$ on $C^\infty(M)$ (see \cite{ABC} and \cite{AC} for further details). Actually, $L_0$ is nothing but the Laplacian on $M$.

As it is shown in \cite[Lemma 3.3]{AC}, if there exists an elliptic point (for a suitable choice of $N$) and the $(k+1)$th mean curvature $H_{k+1}$ satisfies $H_{k+1} > 0$ for some $2 \leq k \leq n-1$, then $L_j$ is elliptic for $1 \leq j \leq k$. Concretely, we can state
\end{rem}

\begin{teor}
\label{teohk}
Let $\psi: M \rightarrow \overline{M}$ be a compact spacelike hypersurface in a GRW spacetime such that $H_{k+1} > 0$ for some $2 \leq k \leq n-1$. If any $H_j$ vanishes nowhere and satisfies $\frac{H_{j+1}}{H_j} \geq \frac{f'(\tau)}{f(\tau)}$ with $f'(\tau)>0$ (resp. $\frac{H_{j+1}}{H_j} \leq \frac{f'(\tau)}{f(\tau)}$ with $f'(\tau)<0$) for $1 \leq j \leq k$, then $M$ is a spacelike slice.
\end{teor}

\begin{demo}
Firstly, we have from \cite[Lemma 5.3]{AC} that there exists an elliptic point of $M$ with respect to an appropiate choice of $N$. Therefore, from this fact and $H_{k+1} > 0$ on $M$ for some $2 \leq k \leq n-1$, \cite[Lemma 3.3]{AC} asserts that the operator $L_i$ is elliptic for all $1 \leq i \leq k$ (using that choice of $N$ if $i$ is odd). Moreover, from 

\begin{eqnarray}
\label{lkf}
L_j(\mathcal{F}(\tau)) &=& -f'(\tau) \mathrm{trace}(P_j) + f(\tau) \overline{g}(N, \partial_t) \mathrm{trace}(A \circ P_j) \nonumber \\
& =& -c_j \big( f'(\tau) H_j + f(\tau) H_{j+1} \overline{g}(N, \partial_t) \big), \nonumber
\end{eqnarray}

\noindent where $c_j \in \mathbb{R}$, \cite[p. 711]{AC}, we get

$$L_j(\mathcal{F}(\tau)) = -c_j f(\tau) H_j \left( \frac{f'(\tau)}{f(\tau)} + \frac{H_{j+1}}{H_j} \overline{g}(N, \partial_t) \right).$$

Hence, taking our assumptions into account we obtain that $L_j(\mathcal{F}(\tau))$ is signed. Since $M$ is compact and $L_j$ is an elliptic operator, the classical maximum principle can be called to conclude that $\mathcal{F}(\tau)$ is constant and $M$ is a spacelike slice.
\end{demo}

\section{Uniqueness results for Einstein GRW spacetimes}
\label{seceins}

We begin this section pointing out that in order to apply the results obtained before to the relevant case of de Sitter spacetime $\overline{M} = \mathbb{R} \times_{\cosh(t)} \mathbb{S}^n$ we need extra assumptions. Hence, directly from Theorem \ref{teoleq} and Proposition \ref{teoHf} we get, respectively

\begin{prop}
\label{prodesit}
The only compact spacelike hypersurfaces with $\tau \leq 0$ (resp., $\tau \geq 0$) in de Sitter spacetime $\overline{M} = \mathbb{R} \times_{\cosh(t)} \mathbb{S}^n$ whose mean curvature function satisfies $H \leq \tanh(\tau)$ (resp., $H \geq \tanh(\tau)$) are the spacelike slices.
\end{prop}

\begin{prop}
\label{promaxdesit}
The only compact spacelike hypersurface with  $\tau \geq 0$ (resp. $\tau \leq 0$) and non-positive (resp., non-negative) mean curvature in de Sitter spacetime $\overline{M} = \mathbb{R} \times_{\cosh(t)} \mathbb{S}^n$ is the totally geodesic spacelike slice $\{0\} \times \mathbb{S}^n$.
\end{prop}

Of course, two nonparametric uniqueness results can be deduced from these propositions where the signs of the involved functions must be suitably assumed. However, as we will show now, if the GRW spacetime is Einstein we can deal with equation (E) without the sign restriction. In order to do so, we will give a new application of an integral formula that was obtained in \cite{A-R-S1}. Namely, for a compact spacelike hypersurface $M$ in a GRW spacetime $\overline{M}$ we have

\begin{equation}
\label{int}
\mathlarger{\int}_M \bigg\{ (n-1) f(\tau) \ \overline{g}(\nabla H, \partial_t) + f(\tau) \ \overline{{\rm Ric}}(\partial_t^T , N) + f(\tau) \ \overline{g}(\partial_t, N) \Big( {\rm trace}(A^2) - n H^2 \Big) \bigg\}dV = 0,
\end{equation}

\noindent where $\nabla H$ represents the gradient of the mean curvature function of $M$. 

Note that a GRW spacetime is Einstein with $\overline{{\rm Ric}} = \overline{c}\thinspace \overline{g}$ if and only if its fiber $(F,g)$ has constant Ricci curvature $c$ and $f$ satisfies the equations

\begin{equation}
\label{eqeins}
\frac{f''}{f} = \frac{c}{n} \hspace{5mm} \text{and} \hspace{5mm} \frac{\overline{c}(n-1)}{n}=\frac{c+(n-1)(f')^2}{f^2}.
\end{equation}

The positive solutions of (\ref{eqeins}) were given in \cite{ARS2}. Furthermore, from (\ref{eqeins}) we also have 

\begin{equation}
\label{logei}
(n-1)(\log f)''= \frac{c}{f^2}.
\end{equation}

Taking into account (\ref{int}) and (\ref{logei}), we can state the following result

\begin{teor}
\label{teointeins}
Let $M$ be a compact spacelike hypersurface in an Einstein GRW spacetime whose mean curvature function verifies $H=\frac{f'(\tau)}{f(\tau)}$. If the constant Ricci curvature of the fiber satisfies $c \geq 0$, then $M$ is totally umbilic. Moreover, if $c > 0$, $M$ must be a spacelike slice.
\end{teor}

\begin{demo}
From our assumptions and (\ref{logei}), the gradient of $H=\frac{f'(\tau)}{f(\tau)}$ in $M$ is

\begin{equation}
\label{gh}
\nabla H = (\log f)''(\tau) \nabla \tau = \frac{c}{(n-1) f^2(\tau)} \nabla \tau.
\end{equation}

Thus, using (\ref{gh}) we can write (\ref{int}) as

\begin{equation}
\label{intein}
-c \int_M \frac{|\nabla \tau|^2}{f(\tau)} dV + \int_M f(\tau) \ \overline{g}(\partial_t, N) \Big({\rm trace}(A^2) - n H^2 \Big) dV = 0,
\end{equation}

\noindent where we have used (\ref{logei}). The result follows using that the Schwarz inequality gives ${\rm trace}(A^2) - n H^2 \geq~0$, with equality holding if and only if $M$ is totally umbilic.
\end{demo}

Particularly, in the case of de Sitter spacetime, the previous result yields to

\begin{coro}
\label{corointds}
The only compact spacelike hypersurfaces in de Sitter spacetime $\overline{M} = \mathbb{R} \times_{\cosh(t)} \mathbb{S}^n$ whose mean curvature function satisfies $H = \tanh(\tau)$ are the spacelike slices.
\end{coro}

Taking into account (\ref{curvmed}) and Corollary \ref{corointds} we are able to get

\begin{coro}
\label{nonpards}
The only entire solutions on $\mathbb{S}^n$ of
$$ {\rm div} \left( \frac{Du}{n \cosh(u) \sqrt{\cosh^2(u) - |Du|^2}} \right) 
 +\frac{\sinh(u)}{n \sqrt{\cosh^2(u)- |Du|^2}} \left( n + \frac{|Du|^2}{\cosh^2(u)} \right) = \tanh(u),$$

$$ |Du| < \cosh(u)$$

\noindent are the constant functions. 

\end{coro}

\section{Uniqueness results for GRW spacetimes obeying NCC}\label{last}
A Lorentzian spacetime obeys the Timelike Convergence Condition (TCC) if its Ricci tensor satisfies
$$ \overline{{\rm Ric}}(X, X) \geq 0,$$
\noindent for all timelike vectors $X$. It is usually argued that TCC is the mathematical way to express that gravity, on average, attracts (see \cite{SW}). Furthermore, if the spacetime satisfies the Einstein equation with a physically reasonable stress-energy tensor, then it must obey TCC \cite[Ex. 4.3.7]{SW}. A weaker energy condition is the Null Convergence Condition (NCC), which reads
$$ \overline{{\rm Ric}}(Z, Z) \geq 0,$$
\noindent for all lightlike vectors $Z$, i.e., $Z \neq 0$ satisfying $\overline{g}(Z,Z)=0$. An easy continuity argument shows that TCC implies NCC. It is well known that a GRW spacetime $\overline{M}=I \times_f F$ satisfies NCC if and only if
\begin{equation}\label{NCCc}
{\rm Ric}^F-(n-1)f^2 (\log f)''\geq 0,
\end{equation}
\noindent where ${\rm Ric}^ F$, denotes the Ricci curvature of the fiber. In particular, if $F$ is compact, then every GRW spacetime with ${\rm Ric}^F>(n-1)\sup\, \big(f^2(\log f)''\big)$ strictly obeys NCC. 

Now, we will make use again of formula (\ref{int}) to state
 
 \begin{teor}
\label{teointeins2}
Let $M$ be a compact spacelike hypersurface in a GRW spacetime whose mean curvature function verifies $H=\frac{f'(\tau)}{f(\tau)}$. If the spacetime obeys NCC and its warping function satisfies $(\log f)''\geq 0$, then $M$ is totally umbilic. Moreover, if NCC is  strictly satisfied or $(\log f)''> 0$ on $M$, then $M$ must be a spacelike slice.
\end{teor}

\begin{demo}
From (\ref{int}), we have

\begin{equation}
\label{iit}
-(n-1)\mathlarger{\int}_M  f(\tau)(\log f)''(\tau)|\nabla \tau|^2 dV + \mathlarger{\int}_M f(\tau) \ \overline{{\rm Ric}}(\partial_t^T , N) dV$$  $$+ \mathlarger{\int}_M f(\tau) \ \overline{g}(\partial_t, N) \Big( {\rm trace}(A^2) - n H^2 \Big) dV = 0.
\end{equation}

Moreover, writing $N = N^F - \overline{g}(\partial_t, N) \partial_t$ and $\partial_t^T = \overline{g}(\partial_t, N) N^F + \left(1-\overline{g}(\partial_t, N)^2 \right) \partial_t$, where $N^F$ denotes the projection of $N$ on the fiber $F$, we get

\begin{equation}
\label{rtn}
\overline{{\rm Ric}}(\partial_t^T , N) = \overline{g}(\partial_t, N) \left\{ \overline{{\rm Ric}}(N^F , N^F) - \left(1-\overline{g}(\partial_t, N)^2 \right) \overline{{\rm Ric}}(\partial_t , \partial_t) \right\}.
\end{equation}

Besides, from \cite[Cor. 7.43]{O'N} we know that

\begin{equation}
\label{rict}
\overline{{\rm Ric}}(\partial_t , \partial_t) = - n \frac{f''(\tau)}{f(\tau)}
\end{equation}

and

\begin{equation}
\label{ricn}
\overline{{\rm Ric}}(N^F , N^F) = {\rm Ric}^F(N^F , N^F) - \left(1-\overline{g}(\partial_t, N)^2 \right) \left( \frac{f''(\tau)}{f(\tau)} + (n-1)\frac{f'(\tau)^2}{f(\tau)^2} \right).
\end{equation}

Hence, using (\ref{rict}) and (\ref{ricn}) in (\ref{rtn}) as well as the fact that $\overline{g}(N^F, N^F) = -1 + \overline{g}(\partial_t, N)^2$, we obtain from (\ref{NCCc}) that the second term in (\ref{iit}) is non-positive under NCC. On the other hand, the first term is non-positive when $(\log f)''\geq 0$. Therefore, $M$ is totally umbilic. Moreover, we also obtain $\overline{{\rm Ric}}(\partial_t^T,N)=0$ on $M$. If NCC is strictly satisfied on $M$ then $\partial_t^T=0$, which means that $M$ is a spacelike slice. The same conclusion is obtained if $(\log f)''> 0$ on $M$.
\end{demo}

\begin{rem}
\label{remphi}
\normalfont
Notice that Theorem \ref{teointeins2} can be extended to a more abstract scenario changing the assumption on the mean curvature of $M$ to $H = \varphi(\tau)$, being $\varphi : I \longrightarrow \mathbb{R}$ smooth and increasing.
\end{rem}

As a consequence of the previous theorem, we obtain the following Calabi-Bernstein type result

\begin{coro}\label{ultimo} Let $F$ be an $n$-dimensional compact Riemannian manifold and let $f:I\longrightarrow ]0,\infty[$ be a smooth function such that $\log f$ is convex. If the Ricci curvature of $F$ is strictly bounded from below by $(n-1)\sup\, \big(f^2(\log f)''\big)$ or $\log f$ is strictly convex, then the only
entire solutions of equation {\rm (E)} are the constant functions.
\end{coro}

\section{Conclusions}

Finally, the results that we have previously obtained enable us to answer our initial question. This is due to the fact that the mean curvature prescription problem that we have solved is equivalent to our initial problem of finding the foliations of a spatially closed GRW spacetime with simply connected fiber by spacelike hypersurfaces whose volume variation is the same that the one measured by the comoving observers for their restspaces.

Therefore, our uniqueness results for equation {\rm (E)} lead to characterize the foliation by spacelike slices by their volume variation in a wide variety of models with physical interest, including de Sitter spacetime (Corollary \ref{nonpards}) as well as other spatially closed GRW spacetime with simply connected fiber where the  warping function is monotone (Corollary \ref{calbeleq}) or where NCC is strictly satisfied (Corollary \ref{ultimo}).

\section*{Acknowledgements}

The authors are partially supported by Spanish MINECO and ERDF project MTM2013-47828-C2-1-P.


\begin{thebibliography}{99}

\bibitem{AA} A.L. Albujer and L.J. Al\'ias, Spacelike hypersurfaces with constant mean curvature in the steady state spacetime, \emph{P. Am. Math. Soc.}, \textbf{137} (2008), 711--721.

\bibitem{ABC} L.J. Al\'ias, A. Brasil Jr. and A.G. Colares, Integral formulae for spacelike hypersurfaces in conformally stationary spacetimes and applications, \emph{P. Edinburgh Math. Soc.}, \textbf{46} (2003), 465--488.

\bibitem{AC} L.J. Al\'ias and A.G. Colares, Uniqueness of spacelike hypersurfaces with constant higher order mean curvature in generalized Robertson-Walker spacetimes, \emph{Math. Proc. Cambridge}, \textbf{143} (2007), 703--729.

\bibitem{A-R-S1} L.J. Al\'ias, A. Romero and M. S\'anchez,
Uniqueness of complete spacelike hypersurfaces of constant mean
curvature in Generalized Robertson-Walker spacetimes, \emph{Gen.
Relat. Gravit.}, \textbf{27} (1995), 71--84.

\bibitem{ARS2} L.J. Al\'ias, A. Romero and M. S\'anchez, Spacelike hypersurfaces of constant mean curvature and Calabi-Bernstein type problems, \emph{T\^ohoku Math. J.}, \textbf{49} (1997), 337--345.

\bibitem{AMM} J.M. Arms, J.E. Marsden and V. Moncrief, The structure of the space of solutions of Einstein's equations. II. Several Killing fields and the Einstein-Yang-Mills equations, \emph{Ann. Phys.}, \textbf{144} (1982), 81--106. 

\bibitem{Azz} A. Azzollini, Ground state solution for a problem with mean curvature operator in Minkowski space, \emph{J. Funct. Anal.}, \textbf{266} (2014), 2086--2095.

\bibitem{BS} R. Bartnik and L. Simon, Spacelike hypersurfaces with prescribed boundary values and mean curvature, \emph{Commun. Math. Phys.} \textbf{87} (1982), 131--152.


\bibitem{BJM} C. Bereanu, P. Jebelean and J. Mawhin, Radial solutions for some nonlinear problems involving mean curvature operators in Euclidean and Minkowski spaces, \emph{P. Am. Math. Soc.}, \textbf{137} (2009), 161--169.

\bibitem{BC} A. Brasil and A.G. Colares, On constant mean curvature spacelike hypersurfaces in Lorentz manifolds, \emph{Mat. Contemp.}, \textbf{17} (1999), 99--136.

\bibitem{BF} D. Brill and F. Flaherty,
Isolated maximal surfaces in spacetime, \emph{Commun. Math. Phys.}
\textbf{50} (1984), 157--165.

\bibitem{CRR} M. Caballero, A. Romero and R.M. Rubio, New Calabi-Bernstein results for some elliptic nonlinear equations, \emph{Anal. Appl.}, \textbf{11} (2013), 1350002(1--13).

\bibitem{Ca} E. Calabi, Examples of Bernstein problems for some nonlinear
equations, \emph{P. Symp. Pure Math.}, \textbf{15} (1970),
223--230.

\bibitem{CY} S.Y. Cheng and S.T. Yau, Maximal space-like hypersurfaces in the
Lorentz-Minkowski spaces, \emph{Ann. of Math.},  \textbf{104} (1976),
407--419.

\bibitem{Ch} Y. Choquet-Bruhat,
Quelques propri\'et\'es des sousvari\'et\'es maximales d'une vari\'et\'e
lorentzienne, \emph{C. R. Acad. Sci. Paris S\'erie A}, \textbf{281} (1975), 577--580.

\bibitem{CCR} C. Corsato, F. Obersnel and P. Omari, Positive solutions of the Dirichlet problem for the prescribed mean curvature equation in Minkowski space, \emph{J. Math. Anal. Appl.} \textbf{405} (2013), 227--239.

\bibitem{dRT} D. de la Fuente, A. Romero and P.J. Torres, Radial solutions of the Dirichlet problem for the prescribed mean curvature equation in a Robertson-Walker spacetime, \emph{Adv. Nonlinear Stud.} \textbf{15} (2015), 171--181.



\bibitem{EH} K. Ecker and G. Huisken, Parabolic methods for the construction of spacelike slices of prescribed mean curvature in cosmological spacetimes, \emph{Commun. Math. Phys.} \textbf{135} (1991), 595--613.

\bibitem{Ge1} C. Gerhardt, H-surfaces in Lorentzian manifolds, \emph{Commun. Math. Phys.} \textbf{89} (1983), 523--553.

\bibitem{Ge2} C. Gerhardt, Hypersurfaces of prescribed mean curvature in Lorentzian manifolds, \emph{Math. Z.} \textbf{235} (2000), 83--97.

\bibitem{JKG} J.L. Jaramillo, J.A.V. Kroon and E. Gourgoulhon, From geometry to numerics: interdisciplinary aspects in mathematical and numerical relativity, \emph{Classical Quant. Grav.}, \textbf{25} (2008), 093001(1--66).

\bibitem{Kaz} J. Kazdan, \emph{Some applications of partial differential equations to problems in Geometry}, Surveys in Geometry Series, Tokyo Univ., (1983).

\bibitem{MT} J.E. Marsden and F.J. Tipler, Maximal hypersurfaces and foliations of constant mean curvature in General Relativity, \emph{Phys. Rep.}, \textbf{66} (1980), 109--139.

\bibitem{O'N} B. O'Neill, \emph{Semi-Riemannian Geometry with
applications to Relativity}, Academic Press, New York, (1983).

\bibitem{PRR1} J.A.S. Pelegr\'in, A. Romero and R.M. Rubio, On maximal hypersurfaces in Lorentz manifolds admitting a parallel lightlike vector field, \emph{Classical Quant. Grav.}, \textbf{33} (2016), 055003(1--8).

\bibitem{PRR2} J.A.S. Pelegr\'in, A. Romero and R.M. Rubio, Uniqueness of complete maximal hypersurfaces in spatially open $(n+1)$-dimensional Robertson-Walker spacetimes with flat fiber, \emph{Gen. Relat. Gravit.}, \textbf{48} (2016), 1--14.

\bibitem{RR} A. Romero and R.M. Rubio, A nonlinear inequality arising in geometry and Calabi-Bernstein type problems, \emph{J. Inequal. Appl.}, (2010), 1--10.

\bibitem{RR2} A. Romero and R.M. Rubio, A nonlinear inequality involving the mean curvature of a spacelike surface in 3-dimensional GRW spacetimes and Calabi-Bernstein type problems, \emph{Contemp. Math.}, \textbf{674} (2016), 141--152.

\bibitem{RRS3} A. Romero, R.M. Rubio and J.J. Salamanca, Spacelike graphs of finite total curvature in certain 3-dimensional generalized Robertson-Walker spacetime, \emph{Rep. Math. Phys.}, \textbf{73} (2014), 241--254.

\bibitem{SW} R.K. Sachs and H. Wu, \emph{General Relativity for Mathematicians}, Graduate Texts in Math., \textbf{48}, Springer-Verlag, New York, (1977).

\bibitem{SY} R. Schoen and S.T. Yau, On the proof of the positive mass conjecture in General Relativity, \emph{Comm. Math. Phys.} \textbf{65} (1979), 45--76.

\bibitem{Se} J. Serrin, The problem of Dirichlet for quasilinear elliptic differential equations with many independent variables, \emph{Philos. T. Roy. Soc. A} \textbf{264} (1969), 413--496.


\end{thebibliography}
\end{document}